\documentstyle[12pt,leqno,amssymb,amscd,graphics]{amsart}

\newtheorem{thm}{Theorem}[section]
\newtheorem{lemma}[thm]{Lemma}

\newtheorem{cor}[thm]{Corollary}
\theoremstyle{definition}
\newtheorem{dfn}[thm]{Definition}
\theoremstyle{remark}

\begin{document}

\newcommand{\ct}{\cite}
\newcommand{\pr}{\protect\ref}
\newcommand{\su}{\subseteq}
\newcommand{\pa}{{\partial}}

\newcommand{\Q}{{\mathbb Q}}
\newcommand{\R}{{\mathbb R}}
\newcommand{\Z}{{\mathbb Z}}
\newcommand{\X}{{\mathbb X}}
\newcommand{\Y}{{\mathbb Y}}

\newcommand{\I}{{\mathrm{Id}}}
\newcommand{\C}{{\mathcal{C}}}

\newcounter{numb}

\title{Complexity of planar and spherical curves}
\author{Tahl Nowik}
\address{Department of Mathematics, Bar-Ilan University, 
Ramat-Gan 52900, Israel}
\email{tahl@@math.biu.ac.il}
\date{February 21, 2008}
\urladdr{www.math.biu.ac.il/$\sim$tahl}

\begin{abstract}
We show that the maximal number of singular moves required to pass between any two 
regularly homotopic planar or spherical curves with at most $n$ crossings, grows quadratically 
with respect to $n$.
Furthermore, this can be done with all curves along the way having at most $n+2$ crossings. 
\end{abstract}

\maketitle

\section{Introduction}\label{intro}

Our subject of interest will be planar and spherical curves,
by which we mean immersions of $S^1$ into $\R^2$ and $S^2$. 
We will be interested in the number of singular occurrences that may be required for
passing between any two regularly homotopic curves, and in the number of crossings
the curves along the way may be required to have.
More precisely, define the distance between two regularly homotopic curves to be the minimal number of singular
occurrences required to pass between them. Define the diameter of a set of regularly homotopic curves to
be the maximal distance between its members. 
We will show that the diameter of the set of all curves in a given regular homotopy class
which have at most $n$ crossings, 
grows quadratically with respect to $n$.
A quadratic upper bound will be established by presenting an explicit algorithm transforming any curve to any 
other regularly homotopic curve. A quadratic lower bound will be established
using the invariant of curves introduced in \ct{n1}.
Furthermore, our explicit algorithm for passing  between curves, has the property that if the initial
and final curves have at most $n$ crossings, then all curves along the way have at most $n+2$ crossings.
It also implies a characterization  of curves with minimal number of crossings in their regular 
homotopy class.

The structure of the paper is as follows. In Section \pr{state} we present the necessary definitions and state our
results. In Section \pr{up} we present our explicit algorithm which proves the quadratic upper bound.
In Section \pr{low} we prove the quadratic lower bound. In Section \pr{con} we compare our results to 
analogous results in other settings, and give a simplified version of our algorithm which may serve as a 
simple constructive proof of Whitney's Theorem.

\section{Definitions and statement of results}\label{state}

A planar or spherical curve is an immersion of $S^1$ into $\R^2$ or $S^2$ respectively.
The Whitney winding number of a planar curve $c$, takes its values in $\Z$, 
and is defined as follows. The derivative of a curve 
$c:S^1 \to \R^2$ does not vanish by definition, and so defines a map 
$S^1 \to \R^2 - \{ 0 \} \simeq S^1$.
The degree of this map is the Whitney winding number of $c$.
For $\omega \in \Z$ we define the curve $\delta_\omega$ to be the
curve with Whitney number $\omega$ described in Figure \pr{cc1}.

For spherical curve $c$, we define the Whitney number of $c$
as 0 or 1 according to the parity of the Whitney number of the planar curve 
obtained by deleting a point from the complement of $c$ in $S^2$.
The curves $\delta_0,\delta_1$,
thought of as spherical curves, will be our chosen representatives
for spherical Whitney number 0 and 1 respectively. 

It is clear that regularly homotopic curves have the same Whitney number. The converse 
is also true, and will follow from our work below. 
We refer to this fact both in the planar and spherical case as Whitney's Theorem.

\begin{figure}[h]
\scalebox{0.8}{\includegraphics{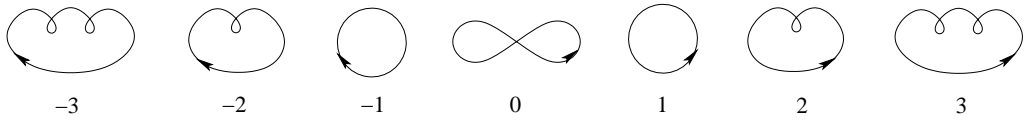}}
\caption{The curves $\delta_\omega$ for $\omega = -3, \dots,3$.}\label{cc1}
\end{figure}

\begin{figure}
\scalebox{0.8}{\includegraphics{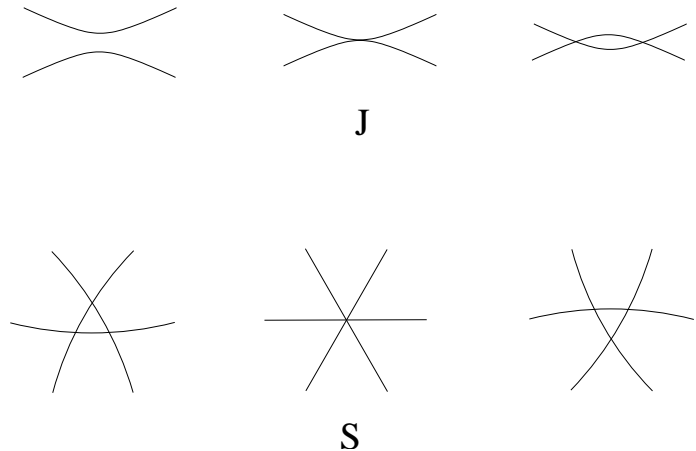}}
\caption{The basic moves.}\label{cc2}
\end{figure}

A curve is called \emph{stable} if its only self intersections are transverse double points.
We denote the space of all stable planar or spherical curves by $\C^P$ and $\C^S$ respectively.
The generic singularities a curve may have are either a tangency of first order between two strands, 
which is called a $J$-type singularity, or three strands meeting at a point, 
each two of which are transverse, which is called an $S$-type singularity.
Singularities of type $J$ and $S$ appear in Figure \pr{cc2}. 
Motion through a $J$ or $S$ singularity as in Figure \pr{cc2} will be called a $J$ or $S$ move,
or jointly, a \emph{basic move}.
An $S$ move preserves the number of crossings of the curve.
On the other hand, the $J$ move which proceeds in Figure \pr{cc2} from left to right 
(respectively from right to left), increases (respectively decreases) the number of crossings by 2.
This will be called an \emph{increasing} $J$ move, and \emph{decreasing} $J$ move, respectively.

\begin{dfn}
\
\begin {enumerate}
\item For two curves $c,c'$ in $\C^P$ or $\C^S$, 
the distance $d(c,c')$ is defined as the minimal number of basic moves
needed to pass from $c$ to $c'$. (If $c,c'$ are not regularly homotopic, set $d(c,c')=\infty$.)
\item For a subset $A$ of $\C^P$ or $\C^S$, the diameter of $A$ is defined as 
$$diam A = \sup_{c,c' \in A} d(c,c').$$
\item Let $B^\omega_n$ be the set of all curves (either in $\C^P$ or $\C^S$) with Whitney number
$\omega$, and with at most $n$ crossings.
\end{enumerate}
\end{dfn}

In this work we prove the following two results.

\begin{thm}\label{t1}
There exist quadratic functions $a(n),b(n)$ (with positive leading coefficients) 
such that for any $\omega$ and $n$,
$$a(n) \leq diam B^\omega_n \leq b(n)$$
\end{thm}

\begin{thm}\label{t2}
For any $c,c' \in B^\omega_n$, there exists a sequence of basic moves from $c$ to $c'$ 
(satisfying the bound of Theorem \pr{t1}), such that all curves along the way are in 
$B^\omega_{n+2}$.
\end{thm}

\section{Upper bound}\label{up}
We present an explicit algorithm for a regular homotopy starting with any $c \in B^\omega_n$ 
and ending with $\delta_\omega$.
In particular this proves Whitney's Theorem.
A simplified version of this algorithm, which may serve as a simple 
constructive proof of Whitney's Theorem, will appear in Section \pr{con}.

\begin{figure}
\scalebox{0.8}{\includegraphics{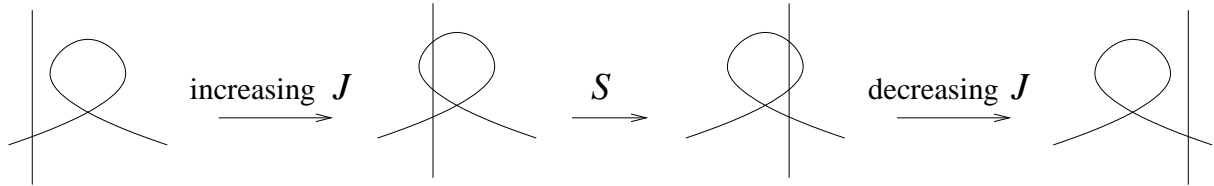}}
\caption{The $Z$ move.}\label{cc3}
\end{figure}

We  define the following \emph{composite} move (as opposed to the \emph{basic} moves $J$ and $S$), 
which we call a $Z$ move, appearing in Figure \pr{cc3}.
This is a sequence of three basic moves: an increasing $J$ move, an $S$ move, and a decreasing
$J$ move. So, during a $Z$ move, the number of crossings increases by 2, but finally decreases
back to the initial number.
Our algorithm will only perform the following moves:
\begin{enumerate}
\item $S$ moves. 
\item Decreasing $J$ moves. 
\item $Z$ moves. 
\end{enumerate}
It follows that the number of crossings of the curves along the way will never exceed $n+2$,
proving Theorem \pr{t2}. 
We present the algorithm for planar curves, and will then present the slight modification required for
spherical curves.

\begin{figure}[h]
\scalebox{0.8}{\includegraphics{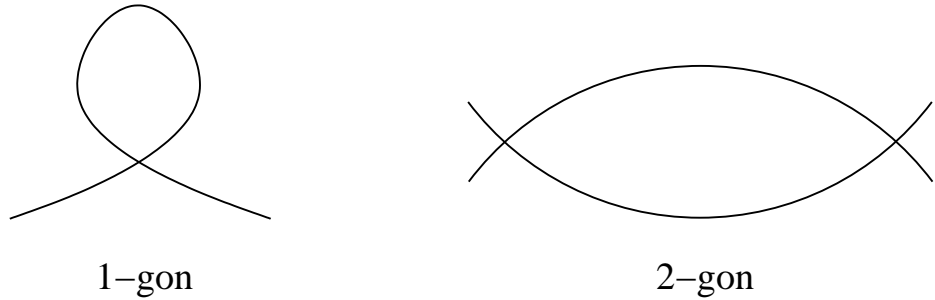}}
\caption{1-gon and 2-gon.}\label{cc4}
\end{figure}

Let $c$ be a \emph{planar} curve. For $k=1,2$, a $k$-gon in $c$, 
is a portion of $c$ of the form appearing in Figure \pr{cc4}. 
The bounded region bounded by the $k$-gon is called its interior. 
The interior of a $k$-gon may intersect other portions of $c$.
It is however assumed as part of the definition, that the corners of a $k$-gon are \emph{convex},
as appears in the figure. This means that the immediate continuations of the arcs 
of the $k$-gon, lie outside its interior.
We will say a $k$-gon in $c$ is \emph{empty}, if its interior is disjoint from $c$.

The algorithm is divided into two steps. Step 1 transforms the initial
curve $c$ into a curve $c'$ which is a string of empty 1-gons. Step 2
transforms $c'$ into $\delta_\omega$. 
The algorithm begins by defining a horizontal line $L$ which crosses $c$
at two points and such that the portion of $c$ below $L$ is an embedded arc.
This is obtained by starting with a horizontal line which is completely
below $c$, and pushing it upwards until slightly after it first touches $c$.
The part of $c$ above and below $L$ will be called the upper curve and lower curve respectively.
Step 1 repeatedly uses two procedures, which we call 
Procedure A and Procedure B. 
Each application of Procedure A or B
reduces the number of crossings of 
the upper curve by at least one, and Procedure A also adds 
an empty 1-gon at one of the ends of the lower curve.
We repeat Procedures A and B 
until there are no crossings in the upper curve, and so $c$ becomes a string of
empty 1-gons.

We note the following:

\begin{lemma}\label{1g}
If the upper curve has crossings, then it contains a 1-gon.
\end{lemma}

\begin{pf}
Parameterizing the upper curve by $[0,1]$, let $t \in [0,1]$ be minimal such
that $c|_{[0,t]}$ is non-injective (i.e. $t$ is the first time that $c$ crosses itself),
and let $s<t$ be such that $c(s)=c(t)$. Then we claim $c|_{[s,t]}$ is a 1-gon.
By definition  of $t$, $c(s)=c(t)$ is the only crossing, and we need to show that the corner at 
$c(s)=c(t)$ is convex.  Let $U \su \R^2$ be the interior of the loop $c([s,t])$.
If the corner is not convex, then  $c(s-\epsilon) \in U$ for some small $\epsilon > 0$,
and so since $c(0) \not\in U$, we must have some $0<r<s-\epsilon$ and $s < r' < t$ with $c(r)=c(r')$,
contradicting the definition of $t$.
\end{pf}

If the upper curve contains an \emph{empty} 1-gon, then we apply 
Procedure A, which is the following.
We slide the 1-gon along the upper curve, 
passing all crossings on the way via $Z$ moves, until it passes $L$ and joins the
lower curve. 
There are two ways that this can be done; the empty 1-gon can be pushed either to the right or to the left.
Let $m$ denote the number of crossings in the upper curve when starting the procedure.  
Since there are $m-1$ crossings that may be passed, and each may be passed twice, 
the sum of the number of $Z$ moves when the 1-gon is pushed to the right and to the left is $2(m-1)$, 
and so either to the right or to the left we need no more than $m-1$ $Z$ moves, which means no more than $3(m-1)$
basic moves. This completes Procedure A.
It reduces the number of crossings in the upper curve by one, 
and adds an empty 1-gon to the lower curve.

If the upper curve contains no empty 1-gons, we apply Procedure B which is the following. 
By Lemma \pr{1g} the upper curve does contain 1-gons. 
Take a minimal such 1-gon, that is, one whose interior $U$ contains no
other 1-gon. By our assumption, this 1-gon is non-empty, and so $U$ contains subarcs of $c$.
If such an arc were not embedded, then by the same argument as in Lemma \pr{1g}, there would be a 1-gon
in $U$, so all such arcs are embedded. It follows that $U$ contains a 2-gon, 
and by taking a minimal one in $U$, we know we have a 2-gon $G$ whose interior $V$
contains no 1-gon and no 2-gon. It follows that $G$ satisfies the following properties.

\begin{enumerate}
\item All arcs in $V$ are embedded.
\item All arcs in $V$ pass from one edge of the 2-gon to the other edge.
\item Each two arcs in $V$ intersect at most once.
\end{enumerate}

We may think of $V$ with its boundary as the square
$[0,1] \times [0,1]$ with the left and right edges $\{ 0 \}\times [0,1]$ and $\{ 1 \}\times [0,1]$
each identified to a point.
And so $V$ itself is parameterized by $(0,1) \times (0,1)$ and so we have a notion of \emph{horizontal
levels} in $V$, namely, the lines $(0,1) \times \{ t \}$.
An arc in $V$ will be called \emph{monotonic} if it passes each horizontal level once 
(as in the definition of braids).  We prove the following:

\begin{figure}
\scalebox{0.8}{\includegraphics{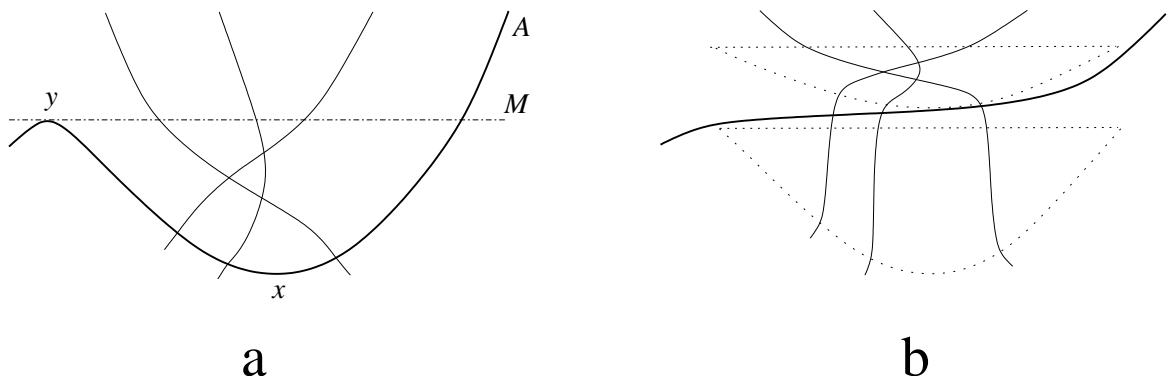}}
\caption{Canceling extrema.}\label{cc5}
\end{figure}

\begin{lemma}\label{gv} 
If $G$ is a 2-gon in $c$ with interior $V$, satisfying properties 1,2,3 above, then 
there is an isotopy of $V$ (keeping a neighborhood of the boundary fixed) after which
all arcs in $V$ become monotonic.
\end{lemma}

\begin{pf}
By induction on the number of arcs. Assume we have $k$ arcs and we have performed
an isotopy making the first $k-1$ arcs monotonic. Denote the $k$th arc by $A$. 
If $A$ is not monotonic, then there are
some extrema along $A$. 
Let $x \in A$ be the highest minimum (by Property 2 the extrema cannot be all maxima). 
If we move from $x$ along $A$ to the right and to the left, then by Property 2 we must reach
some maximum.
If we reach such maximum only on one side (and reach the top edge on the other side)
then let $y$ be that maximum. If we reach a maximum on both sides,
let $y$ be the lower one of these two maxima. Let $M$ be the horizontal line tangent to $A$ at $y$.
It bounds together with the two subarcs ascending from $x$, a sub-region $W$ in $V$. The subarc of $A$ connecting
$x$ and $y$ may not reach $x$ and $y$ from the same side, since if say it reaches both 
$x$ and $y$ from the left, then moving
from $y$ slightly to the right will bring us into $W$, and in order to leave $W$, an additional minimum is required
which is necessarily higher than $x$. So assume when traveling from $x$ to the left, we arrive at $y$ from the right. 
By a level preserving isotopy we may assume the configuration is as in Figure \pr{cc5}a. 
Indeed by Property 3, the upper end of each arc in $W$ must be in $M$.
(The portion above $M$ appearing in the figure, though drawn wide, represents a thin neighborhood above $M$.)
We may now deform a neighborhood 
of $W$ as in Figure \pr{cc5}b to cancel the minimum $x$ with the maximum $y$, keeping the property that 
the first $k-1$ arcs are monotonic. We repeat this process until all extrema of $A$ are canceled, 
and so $A$ is monotonic as well.
\end{pf}

\begin{figure}
\scalebox{0.8}{\includegraphics{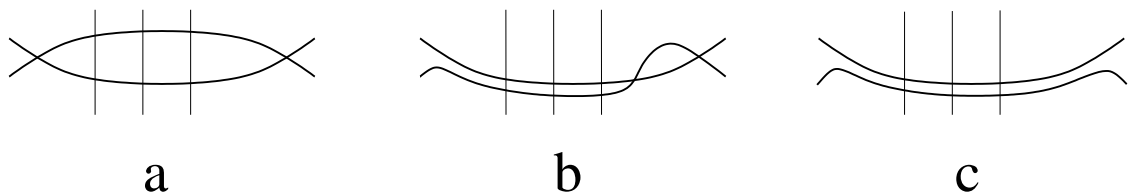}}
\caption{Concluding moves of Procedure B.}\label{cc6}
\end{figure}

So, we now have all arcs in $V$ monotonic, and by an additional slight isotopy we may assume each
crossing in $V$ appears in a different horizontal level. We now move the top edge of $G$ down via the 
horizontal levels, until we pass all crossings. Whenever we pass a crossing, an $S$ move occurs. After this 
sequence of $S$ moves, there are no more crossings in $V$ and so the configuration in $V$ is as in Figure \pr{cc6}a.
Some additional $S$ moves bring us to Figure \pr{cc6}b, and then we perform a single decreasing $J$ move arriving
at Figure \pr{cc6}c, 
by this reducing the number of crossings in the upper curve by 2 (and leaving the lower curve unchanged). 
Let $m$ denote the number of crossings in the upper curve when starting the procedure.
The number of all $S$ moves we have performed is at most $m-2$, and so the total number of 
basic moves is at most $m-1 \leq 3(m-1)$. 
This completes the description of Procedure B. 

As mentioned, we repeat Procedures A and B 
until the upper curve is embedded, and so we have a curve which is a string of empty 1-gons. 
Each performance of Procedure A or B required
at most $3(m-1)$ basic moves, 
and reduced the number of crossings in the upper curve by at least one. 
And so the number of basic moves required to complete Step 1
of the algorithm is at most $3((n-1) + (n-2) + (n-3) + \cdots ) = \frac{3}{2} (n^2 - n)$.

We then begin Step 2 of the
algorithm. If there are two consecutive 1-gons along our string of empty 1-gons, 
which face opposite sides, as on the left hand side of Figure \pr{cc7}, 
then one $Z$ move and one decreasing $J$ move 
as shown in the figure, cancel the two 1-gons. We call this Procedure C.
We repeat Procedure C 
until all 1-gons are facing the same side. There are at most $\frac{1}{2}n$ such pairs of
1-gons, and each such pair was canceled by one $Z$ move and one $J$ move,
which means four basic moves, 
and so this requires at most $2n$ basic moves.

\begin{figure}
\scalebox{0.8}{\includegraphics{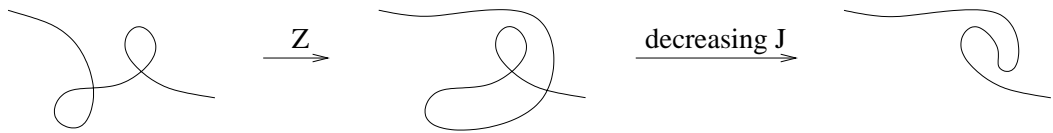}}
\caption{Procedure C: Canceling two consecutive empty 1-gons, facing opposite sides.}\label{cc7}
\end{figure}

\begin{figure}
\scalebox{0.8}{\includegraphics{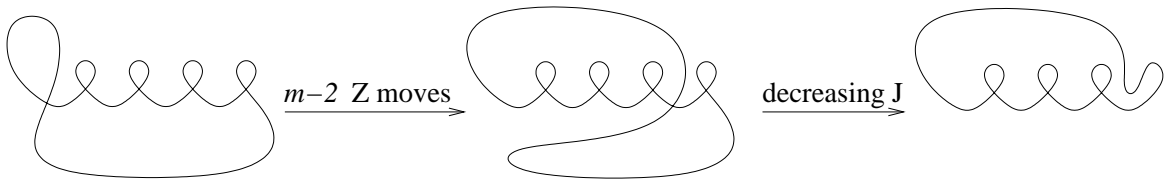}}
\caption{Procedure D: From all 1-gons facing outward, to all 1-gons facing inward.}\label{cc8}
\end{figure}

If at this point all 1-gons are facing inward, or there is just one crossing, then we have reached our 
base curve $\delta_\omega$. Otherwise all 1-gons are facing outward, and there are at least two of them. Let $m$ denote the number of crossings at this time. 
In Figure \pr{cc8} we describe Procedure D, which is 
a sequence of $m-2$ $Z$ moves and one decreasing $J$ move which transforms our curve
into one with all 1-gons facing inward. So this requires $3m-5 \leq 3n-5$ basic moves.

To sum up, our algorithm required at most $\frac{3}{2}n^2 + \frac{7}{2}n - 5$ basic moves. 
It follows that the diameter of $B^\omega_n$ for planar curves is at most $3n^2 + 7n -10$. 

For spherical curves, we perform the exact same algorithm as above until 
obtaining a curve which is a string
of empty 1-gons with all 1-gons facing the same side. 
Since we are now in $S^2$ we can always think of the
1-gons as facing ``outward'', and so if there are at least two 1-gons we can apply
Procedure D which 
requires $3m-5 \leq 3(m-1)$ basic moves and reduces the number of crossings by 2. 
We repeat Procedure D until we have just one or zero
crossings, and so have reached $\delta_0$ or $\delta_1$. This requires at most 
$3((n-1) + (n-3) + (n-5) + \cdots) \leq \frac{3}{4}n^2$ basic moves. 
So all together in the case of spherical curves we need
at most $\frac{9}{4}n^2 + \frac{1}{2}n$ basic moves, 
and so the diameter of $B^\omega_n$ is at most $\frac{9}{2}n^2 + n$.

Our algorithm also provides the following:

\begin{cor}\label{cor}
$\delta_\omega$ has the minimal number of crossings in its regular homotopy class, this minimum being $||\omega|-1|$.
Furthermore, any other curve with the minimal number of crossings can be obtained from $\delta_\omega$ by a sequence of $Z$ moves.
\end{cor}

\begin{pf}
Starting with any curve $c$ of Whitney number $\omega$, our algorithm uses
only $S$ moves, decreasing $J$ moves, and $Z$ moves, and brings us to $\delta_\omega$, which has
$||\omega|-1|$ crossings. Since these three moves do not increase the number of crossings, the first 
statement follows. 
For the second statement, note that of our four procedures, only Procedure A does not decrease the
total number of crossings. 
So if $c$ also has the minimal number of crossings, then our
algorithm will necessarily only apply Procedure A, which only uses $Z$ moves.
\end{pf}

\section{Lower bound}\label{low}

Let $\X$ be the free abelian group with basis all symbols of the form $X_{a,b}$ with $a,b \in \Z$.
In \ct{n1} the invariant $f^X : \C^S \to \X$ is defined as follows.

\begin{dfn}\label{di}
Let $D$ be an oriented 2-disc, and let $e$ be a stable arc in $D$.
\begin{enumerate}
\item For a double point $v$ of $e$
we define $i(v) \in \{1,-1\}$, where $i(v)=1$ if the orientation at $v$ given by the two tangents to $e$ at $v$,
in the order they are visited, coincides with the orientation of $D$. Otherwise $i(v)=-1$.
\item We define the index of $e$, $i(e)\in\Z$ by $i(e) = \sum_v i(v)$ 
where the sum is over all double points $v$ of $e$.
\end{enumerate}
\end{dfn}

For $c \in \C^S$, let $v$ be a double point of $c$, and let $u_1, u_2$ be the two tangents at  
$v$ ordered by the orientation of $S^2$.
Let $U$ be a small neighborhood of $v$ and let $D=S^2 - U$. 
Now $c|_{c^{-1}(D)}$ defines two arcs $c_1,c_2$ in $D$, ordered so that
the tangent $u_i$ leads to $c_i$, $i=1,2$.
We denote $a(v) = i(c_1)$ and $b(v) = i(c_2)$ 
where the orientation on $D$ is that restricted from $S^2$.
Then $f^X : \C^S \to \X$ is defined as follows:
$$f^X(c) = \sum_v X_{a(v),b(v)}$$ 
where the sum is over all double points $v$ of $c$.
We mention that another invariant $f^Y : \C^S \to \Y$ is defined in \ct{n1}, 
and it is shown that the invariant $f^X \oplus f^Y$
is a universal order 1 invariant of spherical curves. 
We define $f^X$ for planar curves in precisely the same way. We mention that it 
may be obtained from the universal order 1 invariant of planar curves appearing in \ct{n2}, 
by the reduction $X^{a_1,b_1}_{a_2,b_2} \mapsto X_{a_1,b_1}$.

Now let $g$ be the $\Z$ valued invariant (for both planar and spherical curves) defined by $g= \phi \circ f^X$
where $\phi : \X \to \Z$ is the homomorphism defined by $X_{a,b} \mapsto a-b$. 
It is shown in \ct{n1} (and applies for planar curves in just the same way)
that the change in the value of $f^X$ due to a basic move is of one of the 
following forms:

\begin{itemize}
\item $X_{a,b} + X_{b,a}$ 
\item $X_{a,b+1} + X_{b,a+1}$
\item $X_{a-1,b}+X_{b-1,a}$
\item $ -X_{a,b+c+2}-X_{b,c+a+2}-X_{c,a+b+2} +X_{a,b+c} +X_{b,c+a} +X_{c,a+b}$
\item $ -X_{c,a+b+1}-X_{b+c-1,a}-X_{b,c+a+1}+X_{b+c+1,a}+X_{b,c+a-1}+X_{c,a+b-1}$
\item $ -X_{c,a+b+1}-X_{c+a-1,b}-X_{b+c-1,a}+X_{b+c+1,a}+X_{c,a+b-1}+X_{c+a+1,b}$
\item $ -X_{b+c-2,a}-X_{c+a-2,b}-X_{a+b-2,c}+X_{a+b,c}+X_{b+c,a}+X_{c+a,b}$
\end{itemize}

Applying $\phi$ to these elements, we see that the change in the value of $g$ due to a basic move is at most 6.
This offers a way for bounding below the distance between two regularly homotopic curves $c,c'$, namely,
$d(c,c') \geq \frac{1}{6}| g(c)-g(c')|$. We use this to obtain a lower bound to the diameter of $B^\omega_n$.
By taking mirror images of all curves, 
it is clear that $diam B^\omega_n = diam B^{-\omega}_n$, and so it is enough to find lower bounds when $\omega \geq 0$.
We note that the number of crossings of a curve is always of opposite parity to its Whitney number.
So if $\omega,n$ are of opposite parity then $B^{\omega}_{n+1} = B^{\omega}_n$ so it is enough to find lower 
bounds for $diam B^{\omega}_n$ for $\omega,n$ of opposite parity.

\begin{figure}[ht]
\scalebox{0.8}{\includegraphics{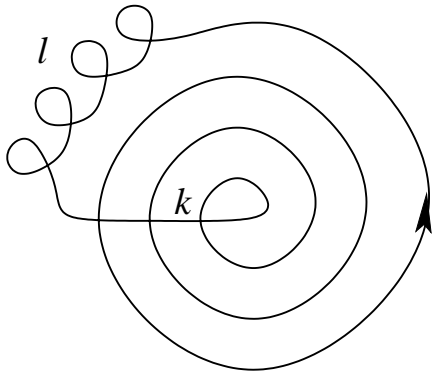}}
\caption{The curve $c_{\omega,n}$ for $\omega = 0$ and $n=7$ ($k=3, l=4$).}\label{cc9}
\end{figure}

For $n \geq \omega \geq 0$ of opposite parity
we construct the curve $c_{\omega,n}$ appearing in Figure \pr{cc9}, having $n$ crossings
and Whitney number $\omega$.
If as indicated in the figure, there are $k$ crossings along the middle horizontal line, 
and $l$ empty 1-gons along the outer circle, then $n=k+l$ and $\omega = k+1-l$.
Inverting these equations gives $k=\frac{1}{2}(n + \omega - 1)$ and $l=\frac{1}{2}(n - \omega +1)$, so 
these are the values for $k$ and $l$ that we use to construct $c_{\omega,n}$. 

Direct inspection of $c_{\omega,n}$ 
gives $$f^X ( c_{\omega,n}) = \sum_{i=1}^k X_{l +i -1 , i-k} + l X_{0,l-1-k}$$
and so $g(c_{\omega,n}) = k(l+k-1) + l(k-l+1) = \frac{1}{2} n^2 + (\omega - 1) n + \frac{1}{2}(1-\omega^2)$.
On the other hand, $f^X(\delta_0) = X_{0,0}$ and for $\omega > 0$, 
$f^X(\delta_\omega) = (\omega - 1)  X_{2-\omega,0}$,
and so $g(\delta_0)=0$ and for $\omega > 0$ $g(\delta_\omega) = -\omega^2 + 3\omega -2$.
Since $diam B^\omega_n = diam B^{-\omega}_n$
we obtain the following lower bound on $diam B^\omega_n$ for any $\omega$, and all $n \geq |\omega|$ of opposite parity, both in $\C^P$ and $\C^S$:
$$diam B^{\omega}_{n+1} =  diam B^\omega_n \geq \frac{1}{12} n^2 + \frac{1}{6}(|\omega| - 1) n + k_\omega \geq \frac{1}{12}(n-1)^2$$  
where $k_0 = \frac{1}{12}$ and for $\omega \neq 0$, $k_\omega = \frac{1}{12}(\omega^2 - 6|\omega| + 5)$,
and the last inequality uses $n \geq |\omega|$.

\section{Concluding remarks}\label{con}

We first sum up our proof of Theorems \pr{t1}, \pr{t2}. The upper and lower bounds given in Sections \pr{up} and \pr{low} establish Theorem \pr{t1}. The fact that our algorithm  
uses only $S$ moves, decreasing $J$ moves, and $Z$ moves, implies Theorem \pr{t2}, since for any 
$c,c' \in B^\omega_n$ use the algorithm to get from $c$ and from $c'$ to $\delta_\omega$.

Next, we remark that the statement of Theorem \pr{t2} cannot be improved. Indeed, it is clear that e.g. our curves 
$c_{\omega,n}$ do not contain a configuration that allows either an $S$ move or a decreasing $J$ move,
and so any regular homotopy from $c_{\omega,n}$ to $\delta_\omega$ must begin with an increasing $J$ move, and
so the number of crossings will increase by 2.

In \ct{my} and \ct{v} a different complexity count is carried out for curves.
The curves are polygonal, and the moves that are counted are translations of a vertex, together with its
two neighboring edges.  It is shown in \ct{v} that the number of such moves required to pass 
between any two regularly homotopic curves (in an appropriate PL sense) 
is linear with respect to the number of edges of the curves.

We would also like to compare our result to what is known for knot diagrams.
One defines the distance between two knot diagrams in terms of the number of Reidemeister moves required
for passing between them. In \ct{hl}, 
an upper bound for the diameter of the set of all diagrams of the unknot with at most $n$ crossings is given
which is exponential with respect to $n$. In \ct{hn}, a lower bound is given which is quadratic with respect to $n$.

Finally, we would like to present a simplified version of our algorithm, which cannot prove 
Theorems \pr{t1}, \pr{t2}, but may serve as a simple constructive proof of Whitney's Theorem,
stating that two curves with the same Whitney winding number are regularly homotopic.
Begin with the same horizontal line $L$, 
and whenever there is an empty 1-gon in the upper curve, proceed as in 
Section \pr{up}.
If there are no empty 1-gon, then by Lemma \pr{1g} there exists a non-empty 1-gon. Shrink this
1-gon towards its corner, until it is empty.
Once there are no crossings in the upper curve, proceed as in Section \pr{up}.
So, this version of the algorithm avoids the complicated procedure B, 
and instead includes the simple procedure of shrinking a non-empty 1-gon.
The bounds on the number of moves and crossings are lost however.

\end{document}